\documentclass[11pt]{article}
\usepackage{amsfonts}
\usepackage{latexsym}
\usepackage{amssymb}
\newcommand{\be}{\begin{equation}} \newcommand{\ee}{\end{equation}}
\newcommand{\bea}{\begin{eqnarray}} \newcommand{\eea}{\end{eqnarray}}
\newcommand{\bean}{\begin{eqnarray*}}
  \newcommand{\eean}{\end{eqnarray*}}
\newcommand{\brray}{\begin{array}} \newcommand{\erray}{\end{array}}
\newcommand{\ben}{\begin{equation}{nonumber}}
  \newcommand{\een}{\end{equation}{nonumber}}
\newcommand{\newsection}[1]{\setcounter{equation}{0}
  \setcounter{dfn}{0}
\section{#1}}

\newtheorem{dfn}{Definition}[section]
\newtheorem{theorem}[dfn]{Theorem}
\newtheorem{lmma}[dfn]{Lemma}
\newtheorem{ppsn}[dfn]{Proposition}
\newtheorem{crlre}[dfn]{Corollary}
\newtheorem{xmpl}[dfn]{Example}
\newtheorem{rmrk}[dfn]{Remark}

\newcommand{\bdfn}{\begin{dfn}}
  \newcommand{\bthm}{\begin{theorem}}
    \newcommand{\blmma}{\begin{lmma}}
      \newcommand{\bppsn}{\begin{ppsn}}
        \newcommand{\bcrlre}{\begin{crlre}}
          \newcommand{\bxmpl}{\begin{xmpl}}
            \newcommand{\brmrk}{\begin{rmrk}}

              \newcommand{\edfn}{\end{dfn}}
            \newcommand{\ethm}{\end{theorem}}
          \newcommand{\elmma}{\end{lmma}}
        \newcommand{\eppsn}{\end{ppsn}}
      \newcommand{\ecrlre}{\end{crlre}}
    \newcommand{\exmpl}{\end{xmpl}} \newcommand{\ermrk}{\end{rmrk}}

 \newcommand{\IZ}{{\mathbb Z}}

\newcommand{\cla}{{\mathcal A}} \newcommand{\clb}{{\mathcal B}}

 \newcommand{\clh}{{\mathcal H}}
 
\newcommand{\clk}{{\mathcal K}}

  \def 
\bbE
{\mbox{\boldmath $E$}}                     

  \def 
\bbe
{\mbox{\boldmath $e$}}

\def\wP{\widetilde{P}} \def\wQ{\widetilde{Q}}

 \def\a*{{\mathcal A}_{h,*}} \def\B{{\mathcal B}(h)}
\def\B1{{\mathcal B}_1(h)} \def\b{{\mathcal B}^{s. a. }(h)} \def\b1{{\mathcal
    B}^{s. a. }_1(h)}

\newcommand{\raro}{\rightarrow}

\begin{document}

\title{C*-algebra generated by projections}
\author{{\sc Partha Sarathi Chakraborty}\\
Indian Statistical Institute\\
7 S.J.S.S Marg, New Delhi  110016, India\\
email: parthac@isid.ac.in\\
partha\_sarathi\_c@hotmail.com}
\date{}
\maketitle
\begin{abstract}
The universal C*-algebra generated by $n$ projections has been described. As an immediate corollary one obtains structure theorem for a pair of projections and the solution to an associated index problem. This puts the study of a pair of projections in the proper perspective of noncommutative geometry.
\end{abstract}
keyword: Index problem,  crossed product, Fredholm module, pair of projections.\\
MSC 46L87.
\section {Introduction}
There has been  considerable amount of work in the operator theory/operator algebra literature  on the structure of a pair of projections. Study of this issue goes back to Dixmier (\cite{DIX}). In  \cite{HAL} Halmos
gives a structure theorem for a pair of projections, which has been rederived  on various occasions. There is also a survey article by C. Davis (\cite{DAV}) on this topic. Avron et. al. \cite{SIM} discusses an index problem in this context. All these approaches use geometric or operator theoretic methods. Our aim is to put these in the proper  perspective of noncommutative geometry (NCG).  To that end we have described the universal C*-algebra generated by $n$ projections as a crossed product. An application of the representation of imprimitivity system yields the  structure theorem for a pair of projections. We believe, our proof of this fact is the shortest and the first one to connect this issue clearly with imprimitivity. Incidentally it may be remarked that Dixmier's paper was reviewed by Mackey and it seems that the aforesaid connection went unnoticed. We also show that an  index problem, the main result  discussed in \cite{SIM}, can be settled easily using general machinery from NCG.
\newsection{Results}
We will refer to unitary operators whose  square is  identity as symmetry.
\bthm
Let $\cla_n$be the universal C*-algebra generated by $n$  projections $P_1, P_2, \cdots ,P_n$. Then we have the following isomorphisms,
$$\cla_n \cong {C^* (\IZ_2)}^{\star n} \cong C^*( F_{n-1}) {\rtimes}  \IZ_2$$
where ${C^* (\IZ_2)}^{\star n}$ is the free product of n copies of the group C*-algebra of $\IZ_2$ and
$C^*( F_{n-1}) {\rtimes} \IZ_2$ is the crossed  product of the group C*-algebra of free group with $n-1$ generators with respect to the $\IZ_2$ action that sends each generator to its inverse.
\ethm
{\bf Proof:} The first isomorphism follows from the following observations,\\
a) ${C^* (\IZ_2)}^{\star n}$ is the universal C*-algebra generated by $n$ symmetries  $U_1, U_2, \cdots ,U_n$;\\
b) $P \mapsto U_P=2P-I$ gives a bijective correspondence between projections and symmetries.\\
For the second isomorphism note that  $C^*( F_{n-1}) {\rtimes}  \IZ_2$ is the universal C*-algebra generated by a symmetry $V$ and $n-1$ unitaries $W_1,\cdots, W_{n-1}$ satisfying $ V W_i V= W_i^{-1}$ for all $i$.
An invertible homomorphism $\pi :C^*( F_{n-1}) {\rtimes}  \IZ_2 \raro {C^* (\IZ_2)}^{\star n}$ is given by
$$ \pi (V)=U_1, \; \; \; \pi (VW_i)=  U_{i+1} \mbox { for } i=1, \cdots ,n-1. $$ \hfill  $\Box$
\brmrk
\rm  We can identify representations of $\cla_2$ with that of $C(S^1)$ on some Hilbert space $\clh$ along with a symmetry $V$ implementing the conjugation action of $\IZ_2$ on $S^1$. Clearly $\clh$ must look like $\clk_{-}\oplus \clk_{+}\oplus L^2(S_{+},\mu_{+}) \oplus L^2(S_{-},\mu_{-})$ where $S_{+}=\{
z \in S^1 : Im z >0\}$,   $S_{-}=\{
z \in S^1 : Im z <0\}$ and the measures $\mu_{\pm}$ satisfy $\mu_{+}=\mu_{-} \circ Conjugation$. The symmetry $V$ switches $L^2(S_{\pm},\mu_{\pm})$ and respects $\clk_{\pm}$. The representation $\pi : C(S^1) \raro \clh$ acts as $\pi(f)=$ multiplication by $f$  on  $L^2(S_{\pm},\mu_{\pm})$ and
 $\pi(f)|_{\clk_{\pm}}=f(\pm1)I_{\clk_{\pm}}$. If we denote by $\wP_1,\wP_2$ the images of the generators $P_1,P_2$ of $\cla_2$, then they respect $\clk _{-}\oplus \clk_{+}$ and commute there. Using the correspondence between projections and symmetries stated in the proof of theorem 2.1 we get
that on $L^2(S_{+},\mu_{+}) \oplus L^2(S_{-},\mu_{-})$, $\wP_1= \frac{1}{2}\left( \matrix{ 1 & 1 \cr 1  & 1 \cr } \right),$ and $ \wP_2=\frac{1}{2} \left( \matrix{ 1 & \pi(z)  \cr \pi(z)   & 1 \cr } \right)$. This is the structure theorem for a pair of projections. Note that in this part of the space $\wP_1-\wP_2$ is an offdiagonal operator matrix and on $\clk _{-}\oplus \clk_{+}$ has eigenvalues $\pm1,0$. Therefore under the assumption that
$(\wP_1-\wP_2)^{2k+1}$ is trace class for some $k$, we have the following equality, $tr_{\clh} (\wP_1-\wP_2)^{2k+1+2m}=tr_{\clh}(\wP_1-\wP_2)^{2k+1}, \: \forall m \ge 0 .$
\ermrk
Now we intend to put the index problem tackled by Avron et. al. (\cite{SIM}) in the framework  of NCG. There it has been remarked that the problem has ``relations to the work of Connes \cite{CON1}". We show that the main result of \cite{SIM} can be deduced as a simple application from results in \cite{CON1}.
\bthm
Let $ \wP,\wQ$ be  projections on $\clh$ such that $\wP-\wQ \in \clb_{2k+1}$, the Schatten ideal. Then $\wQ\wP:\wP\clh \raro \wQ\clh$ is Fredholm with $Index (\wQ\wP)= tr (\wP-\wQ)^{2k+1}$.
\ethm
{\bf Proof:} Let $\clk= \clh \oplus \clh$, $\gamma= \left( \matrix{ 1 & 0 \cr 0  & -1 \cr } \right),
F=\left( \matrix{ 0& 1 \cr 1 &  0 \cr } \right)$, and $\pi: \cla_2 \raro \clb(\clk)$ be the representation given on the generators $P_1,P_2$ by the prescription  $\pi(P_1)= \left( \matrix{ \wP & 0 \cr 0  & \wQ \cr } \right), \pi(P_2)= \left( \matrix{ \wQ & 0 \cr 0  & \wP \cr } \right)$.
Then $(\cla_2, \pi, \clk, \gamma, F)$ is an even Fredholm module.  Hence by proposition 2, page 289, and proposition 4, page 296, \cite{CON} we get  $\wQ\wP:\wP\clh \raro \wQ\clh$ is Fredholm with $Index (\wQ\wP)= (-1)^{k+1} tr_\clk \gamma  \pi(P_1) [F, \pi(P_1)]^{2k+2}=
tr (\wP-\wQ)^{2k+3}$. Now  remark 2.2 completes the proof. \hfill $\Box$

  \end{document}